\documentclass[conference]{IEEEtran}
\IEEEoverridecommandlockouts
\usepackage{cite}
\usepackage{amsmath,amssymb,amsfonts}
\usepackage{algorithmic}
\usepackage{graphicx}
\usepackage{textcomp}
\usepackage{xcolor}
\def\BibTeX{{\rm B\kern-.05em{\sc i\kern-.025em b}\kern-.08em
    T\kern-.1667em\lower.7ex\hbox{E}\kern-.125emX}}

\def\bo0{\boldsymbol{0}}

\def\boa{\boldsymbol{a}}

\def\bop{\boldsymbol{p}}
\def\boq{\boldsymbol{q}}
\def\bov{\boldsymbol{v}}
\def\boy{\boldsymbol{y}}
\def\botheta{\boldsymbol{\theta}}

\def\botw{\boldsymbol{\tilde{w}}}

\def\D{^{\rm D}}

\def\G{^{\rm G}}
\def\Gmin{^{\rm G, min}}
\def\Gmax{^{\rm G, max}}

\def\Gprior{^{\rm G, act}}
\def\Gstep{^{\rm G, step}}
\def\GstepOpt{^{\rm G, step, opt}}
\def\Gstepprior{^{\rm G, step, act}}

\def\omin{^{\rm min}} % \min is already defined/assigned!
\def\omax{^{\rm max}} % \max is already defined/assigned!

\def\kK{{\cal K}}
\def\kL{{\cal L}}

\def\kN{{\cal N}}

\def\boQ{\boldsymbol{Q}}

\def\boX{\boldsymbol{X}}

\def\boldx{\boldsymbol{x}} % \box is not possible --> \boldx

\def\boz{\boldsymbol{z}}
\def\yS_kl{y_{{\rm S}, k, l_k}}
\def\yISI_kl{y_{{\rm ISI}, k, l_k}}
\def\yMUI_k{y_{{\rm MUI}, k}}
\def\yN_k{y_{{\rm N}, k}}
\def\P_max{P_{\rm UB}}
\def\Pmin_kl{P_{{\rm min},k,l_k}}
\def\p_i_max{p_{i, \rm max}}
\def\bQ_ISIkl{\boQ_{\rm ISI, k,l_k}}
\def\cQ_MUIk{\boQ_{\rm MUI, k}}
\def\dQ_Nk{\boQ_{\rm N, k}}
\def\z_opt{\boz_{\rm opt}}
\def\tbw_opt{\botw_{\rm opt}}
\def\X_opt{\boX_{\rm opt}}

\def\redtext{}

\begin{document}

\title{\redtext{A Computationally Efficient Method for Solving Mixed-Integer AC Optimal Power Flow Problems}\\

\thanks{
% \noindent
\begin{minipage}{0.2\columnwidth}
\includegraphics[width=\columnwidth]{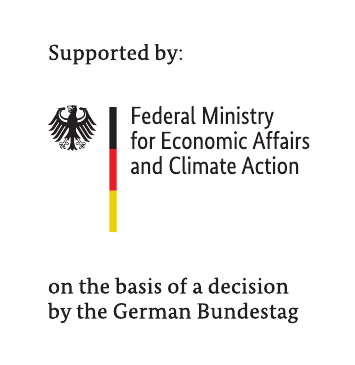}
\end{minipage}
\begin{minipage}{0.75\columnwidth}
This work was supported by the Federal Ministry of Economic Affairs and Climate Action on the basis of a decision by the German Bundestag within the project "Redispatch 3.0" (FKZ 03EI4043J). 
\end{minipage}
}
}

\author{
\IEEEauthorblockN{1\textsuperscript{st} Johannes Heid}
\IEEEauthorblockA{\textit{Sustainable Electrical Energy Systems} \\
\textit{University of Kassel, Fraunhofer IEE}\\
Kassel, Germany \\
johannes.heid@uni-kassel.de}
\and
\IEEEauthorblockN{2\textsuperscript{nd} Nils Bornhorst}
\IEEEauthorblockA{\textit{Sustainable Electrical Energy Systems} \\
\textit{University of Kassel}\\
Kassel, Germany \\
nils.bornhorst@uni-kassel.de}
\and
\IEEEauthorblockN{3\textsuperscript{rd} Eric Tönges}
\IEEEauthorblockA{\textit{Sustainable Electrical Energy Systems} \\
\textit{University of Kassel}\\
Kassel, Germany \\
eric.toenges@uni-kassel.de}
\and
\IEEEauthorblockN{4\textsuperscript{th} Philipp Härtel}
\IEEEauthorblockA{\textit{Sustainable Electrical Energy Systems} \\
\textit{University of Kassel, Fraunhofer IEE}\\
Kassel, Germany \\
philipp.haertel@iee.fraunhofer.de}
\and
\IEEEauthorblockN{5\textsuperscript{th} Denis Mende}
\IEEEauthorblockA{\textit{Sustainable Electrical Energy Systems} \\
\textit{University of Kassel, Fraunhofer IEE}\\
Kassel, Germany \\
denis.mende@iee.fraunhofer.de}
\and
\IEEEauthorblockN{6\textsuperscript{th} Martin Braun}
\IEEEauthorblockA{\textit{Sustainable Electrical Energy Systems} \\
\textit{University of Kassel, Fraunhofer IEE}\\
Kassel, Germany \\
martin.braun@uni-kassel.de}
}

% % replaced by this from the PSCC paper (To Do: Check if this is allowed!):
% \author{\IEEEauthorblockN{Johannes Heid\IEEEauthorrefmark{1},
% Nils Bornhorst\IEEEauthorrefmark{1},
% Eric Tönges\IEEEauthorrefmark{1},
% Philipp Härtel\IEEEauthorrefmark{1}\IEEEauthorrefmark{4}, 
% Denis Mende\IEEEauthorrefmark{1}\IEEEauthorrefmark{4}, and 
% Martin Braun\IEEEauthorrefmark{1}\IEEEauthorrefmark{4}}
% \IEEEauthorblockA{\IEEEauthorrefmark{1}
% University of Kassel, Department e$^2$n,
% Kassel, Germany\\}
% % \IEEEauthorblockA{\IEEEauthorrefmark{2} Department Name of Organization B\\
% % Name of the organization B,
% % Address B\\ Emails if wanted}
% % \IEEEauthorblockA{\IEEEauthorrefmark{3} Department Name of Organization C\\
% % Name of the organization C,
% % Address C\\ Emails if wanted}
% \IEEEauthorblockA{\IEEEauthorrefmark{4}
% Fraunhofer IEE,
% Kassel, Germany\\ johannes.heid@uni-kassel.de}
% }

\maketitle

\begin{abstract}
\redtext{
Stepwise controllable devices, such as switched capacitors or stepwise controllable loads and generators, transform the nonconvex AC optimal power flow (AC-OPF) problem into a nonconvex mixed-integer (MI) programming problem which is generally hard to solve optimally. Existing methods for solving MI-AC-OPF problems usually suffer from either limited accuracy or computational intractability, making them impractical for real-world applications. To address these challenges, we propose an efficient iterative deflation approach providing high-quality approximate solutions. In each iteration, a continuously relaxed version of the MI-AC-OPF problem is solved and one candidate integer value is systematically eliminated based on the evaluation of a simple power flow result. The computational complexity of the proposed algorithm grows linearly with the number of integer optimization variables, ensuring scalability. Simulations demonstrate that the proposed approach achieves significant improvements in solution accuracy compared to a state-of-the-art approach. Thus, the proposed method is promising for solving practical MI-AC-OPF problems.
}
\end{abstract}


\section{Introduction}\label{sec:intro}

% To Do: Check the international definition of DSO! Probably, HV grids are operated by TSOs in the US!
\subsection{Motivation}
In the ongoing energy transition, distribution system operators (DSOs) are facing a rapidly increasing number of distributed generators (DGs) (e.g., PV systems and wind farms) and electrified loads (e.g., heat pumps and electric vehicles) connected to the distribution grid. Since grid expansion cannot keep pace with this rapid development, congestion management has become a crucial issue (see \cite{Gross} and references therein). At the same time, DGs are commonly controllable
\redtext{
and can thus provide flexibilities to enable congestion management in the distribution grid. These flexibilities of the controllable DGs can be used as optimization variables in an AC optimal power flow (AC-OPF) problem to mitigate congestions. However, stepwise controllable devices, such as legacy wind farms with stepwise active power curtailment, switched capacitors, and on-load tap changers, are (still) widely deployed in HV grids \cite{Gross} - \cite{Savasci}. The discrete setpoints of these devices transform the inherently nonconvex AC-OPF problem into a mixed-integer (MI) nonconvex problem which is generally hard to solve optimally.
It is worth noting that a DC approximation of the AC-OPF problem exists, simplifying the problem into a linear DC-OPF formulation. However, the DC approximation is only applicable at the EHV level, where line resistances can be neglected. At the HV level (and the MV and LV levels), the resistive component of the line series impedance cannot be neglected, making the AC-OPF formulation necessary.
}

\subsection{State-of-the-art}

Approaches to solve MI-AC-OPF problems for distribution grids with stepwise controllable devices (here of load tap changers and switched capacitors) has been the subject of recent research \cite{Shukla}, \cite{Savasci}. Solvers for convex MI programming problems, such as MI linear programming (LP), MI second-order cone programming (SOCP), and MI semi-definite programming (SDP) problems, exist \cite{Hijazi}. Applying MI programming solvers to the MI-AC-OPF problem has several drawbacks: An approximation or relaxation of the original MI-AC-OPF problem is required to turn the nonconvex problem into a convex one \cite{Molzahn}. When the approximation or relaxation is supposed to be more accurate than a linearization, \text{MI-SOCP} or \text{MI-SDP} solvers have to be used as in \cite{Shukla} and \cite{Savasci}. These solvers rely on branch-and-bound (BnB) or branch-and-cut (BnC) procedures and the overall computational complexity is still very high, growing exponentially with the number of integer variables. Such a computational burden is usually prohibitively high for real-time applications including curative congestion management \cite{Hijazi}.
The commercial mathematical programing language AMPL \cite{Fourer} offers a broad library of solvers, including solvers such as Knitro which are able to approximately solve the original nonconvex MI-nonlinear programming (MINLP) problem directly. However, to the best of our knowledge, all of AMPL's MI solvers rely on BnB or BnC algorithms.

An alternative, straight-forward approach is a two-step algorithm, where a continuously relaxed version of the MI-AC-OPF problem is approximately solved using existing solvers, such as IPOPT (see \cite{Molzahn} and references therein), the continuously relaxed variables are rounded to the closest discrete setpoint and fixed, and finally, the AC-OPF problem is solved again to re-optimize the continuous variables. The results, however, may severely lack accuracy since only two steps are conducted. The two-step approach has been used before for a joint transformer tap position and DG feed-in optimization in \cite{Marten} and \cite{Mamandur}. 
In \cite{Paudyal}, all possible combinations of rounding down or up to the nearest smaller or larger discrete setpoint are considered. For each candidate combination, the AC-OPF problem is solved again for the continuous variables. The combination yielding the best solution is then selected. Compared to the two-step approach, this approach improves the solution accuracy. However, the number of combinations and, thus, the overall computational burden grows exponentially with the number of discrete variables, limiting its practicability for large-scale problems.

In the context of wireless communications, an iterative deflation approach has been proposed for the joint optimization of continuous filter-and-forward relay filter coefficients and discrete decoding delays \cite{Bornhorst2015}, \cite{Bornhorst2013}. The name \textit{deflation approach} stems from the fact that in each iteration of the algorithm, a single decoding delay is removed from the set of all possible decoding delays, until only one decoding delay remains. After each removal, a continuously relaxed version of the original MI programming problem is solved. This iterative process allows for a suboptimal yet near-optimal solution to the MI programming problem to be approached systematically. The computational complexity of the proposed scheme grows only linearly with the number of integer variables. However, in \cite{Bornhorst2013}, it is demonstrated by simulations that the solution obtained with the deflation approach is close to a lower bound and, consequently, near-optimal.

\subsection{Own contribution}

To improve the results of the two-step approach, we propose an efficient iterative deflation approach inspired by the one in \cite{Bornhorst2015} and \cite{Bornhorst2013}. In the proposed approach, one candidate discrete setpoint is eliminated in each iteration. The approach from \cite{Bornhorst2015} and \cite{Bornhorst2013} has to be severely adapted to the new context of power systems.
% To Do: Keep this?
In \cite{Bornhorst2015} and \cite{Bornhorst2013}, the decision to eliminate a candidate discrete value is based on the calculation of the signal-to-interference-plus-noise ratio (SINR) which simply exists in closed-form. In contrast, a novel approach, where the elimination decisions is based on the evaluation of a simple power flow result, is proposed in this work. Specifically, the candidate yielding the strongest constraint violations (voltage band and line loading) and the worst objective function value is eliminated.
% To Do: Keep this?
We introduce a metric for the evaluation of constraint violations and an approach to combine the evaluation of constraint violations with the objective function value to achieve near-optimal results.

% To Do: Is this complete?
% In this paper, we consider the problem of minimizing the total amount of curtailment of DGs and of control of controllable loads that is required to resolve all congestions in the grid.

% Simulation results demonstrate that the proposed technique outperforms the “two-step” approach
% %and MILP-based approaches
% in terms of accuracy at slightly increased computational complexity. It is hence an appropriate trade-off for practical applications.

The remainder of this paper is organized as follows. After defining the nomenclature in Section~\ref{sec:nomenclature}, the considered optimization problem is formulated in Section~\ref{sec:problem}. The state-of-the-art two-step approach is revisited in Section~\ref{sec:2step} and the proposed deflation approach is introduced in Section~\ref{sec:deflat}. Simulation results are presented in Section~\ref{sec:simulations} and conclusions are drawn in Section~\ref{sec:conclusion}.

\section{Nomenclature}\label{sec:nomenclature}

\subsection{Sets}
\addcontentsline{toc}{section}{Nomenclature}
\begin{IEEEdescription}[\IEEEusemathlabelsep\IEEEsetlabelwidth{$\kL_k = \{1, \ldots, L_k\}$}]
\item[$\kN = \{1, \ldots, N\}$] Set of all $N$ buses in the grid
\item[$\kN_m$] Set of buses directly connected to bus $m$
\item[$\kK = \{1, \ldots, K\}$] Set of all $K$ stepwise controllable generators in the grid
% \item[$\kN\step = \{m_1, \ldots, m_K\}$] Set of all $K$ buses with stepwise controllable generators
\item[$\kK_m$] Set of all stepwise controllable generators at bus $m$
\item[$\kL_k = \{1, \ldots, L_k\}$] Set of all $L_k$ discrete setpoints of the $k$th stepwise controllable generator
% \item[$\kL_m = \{1, \ldots, L_m\}$] Set of all $L_m$ discrete setpoints of the stepwise controllable generator connected to bus $m$
\end{IEEEdescription}

\subsection{Continuous optimization variables}
\addcontentsline{toc}{section}{Nomenclature}
\begin{IEEEdescription}[\IEEEusemathlabelsep\IEEEsetlabelwidth{$q\Gstep_k$}]
\item[$p\G_m, q\G_m$] (Re-)active power of continuously controllable generator at bus $m$
% \item[$q\G_m$] Reactive power of continuously controllable generator at bus $m$
\item[$q\Gstep_k$] Reactive power of $k$th stepwise controllable generator 
\item[$v_m$] Voltage magnitude at bus $m$ 
\item[$\theta_m$] Voltage angle at bus $m$
\end{IEEEdescription}

\subsection{Discrete optimization variables}\label{subsec:discrVars}
\addcontentsline{toc}{section}{Nomenclature}
\begin{IEEEdescription}[\IEEEusemathlabelsep\IEEEsetlabelwidth{$p\Gstep_k \in \{p\Gstep_{k,1}, \ldots, p\Gstep_{k,L_k}\}$}]
\item[$p\Gstep_k \in \{p\Gstep_{k,1}, \ldots, p\Gstep_{k,L_k}\}$] Active power of $k$th stepwise controllable generator
\end{IEEEdescription}

\subsection{Auxillary variables:}
\addcontentsline{toc}{section}{Nomenclature}
\begin{IEEEdescription}[\IEEEusemathlabelsep\IEEEsetlabelwidth{$a_{k,l_k} \in \{0, 1\}$}]
\item[$p_{mn}, q_{mn}$] (Re-)active power flow from bus $m$ to bus $n$
% \item[$q_{mn}$] Reactive power flow from bus $m$ to bus $n$
\item[$s_{mn}$] Apparent power flow from bus $m$ to bus $n$
\item[$a_{k,l_k} \in \{0, 1\}$] Binary variables for reformulating the discrete optimization variables in \ref{subsec:discrVars}.
\end{IEEEdescription}

\subsection{Constants/parameters:}
% \begin{supertabular}{l p{0.65\columnwidth}}
\addcontentsline{toc}{section}{Nomenclature}
\begin{IEEEdescription}[\IEEEusemathlabelsep\IEEEsetlabelwidth{$p\Gmin_m, p\Gmax_m$}]
\item[$p\Gprior_m$] Active power of continuously controllable generator at bus $m$ before the optimization
\item[$p\Gstepprior_k$] Active power of $k$th stepwise controllable generator before the optimization
\item[$p\Gstep_{k,l_k}$] $l_k$th discrete active power setpoint of $k$th stepwise controllable generator
\item[$p\D_m, q\D_m$] (Re-)active power of (uncontrollable) load at bus $m$
% \item[$q\D_m$] Reactive power of (uncontrollable) load at bus $m$
\item[$G_{mn}, B_{mn}$] Real and imaginary part of the $m,n$th entry of the bus admittance matrix
\item[$v\omin_m, v\omax_m$] Minimum and maximum voltage magnitude at bus $m$ (voltage band)
\item[$s\omax_{mn}$] Maximum line loading at line connecting buses $m$ and $n$
\item[$p\Gmin_m, p\Gmax_m$] Minimum and maximum active power of continuously controllable generator at bus $m$
\item[$\alpha$] reactive to active power ratio %cos($\varphi$)}
% \item[$N$] number of buses in the grid
% \item[$L_k$] number of discrete setpoints of the $k$th stepwise controllable device
% \item[$K\G$] number of stepwise controllable generators in the grid
% \item[$K\D$] number of stepwise controllable loads in the grid
% \end{supertabular}
\end{IEEEdescription}

\section{Problem formulation}\label{sec:problem}

\redtext{
We consider congestion management in the HV grid, although the formulations are generally applicable to other voltage levels. The goal is to resolve all congestions while minimizing the total active power curtailment of DGs. In German distribution grids, legacy DGs, particularly wind farms, which are limited to be curtailed to 0\,\%, 30\,\%, and 60\,\% (and 100\,\%) of their nominal power, remain widespread \cite{Gross}. Accordingly, we assume that the active power feed-in of some DGs is continuously controllable while others are limited to discrete curtailment steps. Based on the above considerations, the resulting MI-AC-OPF problem can be formulated as follows.
}

% To Do: Replace eqnarray by IEEEeqnarray (everywhere)
\redtext{
\allowdisplaybreaks
\begin{subequations}\label{eqn:problemOrig}
\begin{eqnarray}
  & & \min_{\boldx, \bop\Gstep, \boq\Gstep} f_1(\bop\G, \bop\Gstep) \label{subeqn:objective} \\
  & & \text{subject to} \nonumber \\  
  & & p\G_m + \sum_{k \in \kK_m}p\Gstep_k - p\D_m = \sum_{n \in \kN_m} p_{mn}, \; \forall m \in \kN \label{subeqn:1stLFconstr} \\
  & & q\G_m + \sum_{k \in \kK_m}q\Gstep_k - q\D_m = \sum_{n \in \kN_m} q_{mn}, \; \forall m \in \kN \label{subeqn:2ndLFconstr} \\
  & & p_{mn} = - v_m v_n (G_{mn}\cos{\theta_{mn}} + B_{mn} \sin{\theta_{mn}}) \nonumber \\
  & & + v_n^2 G_{mn}, \; \forall m, n \in \kN \\
  & & q_{mn} = - v_m v_n (G_{mn}\sin{\theta_{mn}} - B_{mn} \cos{\theta_{mn}}) \nonumber \\
  & & - v_n^2 B_{mn}, \; \forall m, n \in \kN \label{subeqn:lastLFconstr}\\
  & & v\omin_m \leq v_m \leq v\omax_m, \; \forall m \in \kN \label{subeqn:firstIneqConstr} \\
  & & p_{mn}^2 + q_{mn}^2 \leq (s\omax_{mn})^2, \; \forall m, n \in \kN  \\
  & & p\Gmin_m \leq p\G_m \leq p\Gmax_m, \; \forall m \in \kN \label{subeqn:lastLimitConstr} \\
  & & q\G_m = \alpha\G_m \cdot p\G_m, \; \forall m \in \kN \label{subeqn:cosPhiConstr} \\
  & & q\Gstep_k = \alpha\Gstep_k \cdot p\Gstep_k, \; \forall k \in \kK \label{subeqn:lastContConstr} \\
  & & p\Gstep_k \in \{p\Gstep_{k,1}, \ldots, p\Gstep_{k,L_k}\}, \; \forall k \in \kK \label{subeqn:stepConstr1}
\end{eqnarray}
\end{subequations}
where
\begin{subequations}
 \begin{eqnarray}
    \boldx & = & [\boy^T, \bov^T, \botheta^T]^T \\
    \boy & = & [(\bop\G)^T, (\boq\G)^T]^T \\
    \bov & = & [v_1, \ldots, v_N]^T \\
    \botheta & = & [\theta_{11}, \ldots, \theta_{NN}]^T, \; \theta_{mn} = \theta_m - \theta_n \\
    \bop\G & = & [p\G_1, \ldots, p\G_N]^T \\
    \boq\G & = & [q\G_1, \ldots, q\G_N]^T \\
    \bop\Gstep & = & [p\Gstep_1, \ldots, p\Gstep_K]^T \\
    \boq\Gstep & = & [q\Gstep_1, \ldots, q\Gstep_K]^T
\end{eqnarray}   
\end{subequations}
and the objective function $f_1(\cdot)$ is defined as follows
\begin{eqnarray}\label{eqn:objective}
  & & f_1(\bop\G, \bop\Gstep) = \sum_{m \in \kN}|p\Gprior_m - p\G_m| \nonumber \\
  & & + \sum_{k \in \kK}|p\Gstepprior_k - p\Gstep_k|.
\end{eqnarray}
}

\redtext{
Generator curtailment must be reimbursed by the DSO, resulting in additional costs for the DSO. Consequently, the DSO's objective is to minimize these costs. The actual costs are unknown to the authors. However, it is reasonable to assume that they are proportional to the amount of curtailed power. Therefore, we approximate the objective of minimizing the curtailment costs by focusing on minimizing the total amount of curtailed power. This is reflected in the objective function \eqref{eqn:objective} of problem \eqref{eqn:problemOrig} which
aims at minimizing the total deviation of the optimized active power feed-in from the actual ones. Further, \eqref{subeqn:1stLFconstr} - \eqref{subeqn:lastLFconstr} represent the power flow equations, \eqref{subeqn:firstIneqConstr} - \eqref{subeqn:lastLimitConstr} represent the limits for the voltage band, the line loading and the continuously controllable generators, respectively. \eqref{subeqn:cosPhiConstr} and \eqref{subeqn:lastContConstr} force the reactive power generated by both, the continuously and the stepwise controllable generators, to follow a fixed cos($\varphi$). Note that this is a reasonable assumption for stepwise controllable generators since they are able to receive only a limited amount of analogue control signals. E.g., those generators can receive one control signal for each of four discrete active power setpoints. Additional control signals for reactive power are usually not available. For simplicity, we make the same assumption for continuously controllable generators. However, \eqref{subeqn:cosPhiConstr} could straightforwardly be replaced by inequality contraints for continuously controllabele $q\G_m$. Further note that $\alpha$ can straightforwardly be converted to cos($\varphi$). Finally, \eqref{subeqn:stepConstr1} represents the possible discrete setpoints to which the stepwise controllable generators can be curtailed. Note that in \eqref{subeqn:stepConstr1}, the curtailment steps are arranged in ascending order with $p\Gstep_{k, 1}$ being the lowest curtailment step (usually = 0) and $p\Gstep_{k, L_k} = p\Gstepprior_k$, i. e., the DGs can either remain uncurtailed or be curtailed to a discrete setpoint below the actual active power feed-in. In case there is no stepwise controllable DG at bus $m$, $\kK_m = \emptyset$. In case there is no continuously controllable DG at node $m$, $p\Gmin_m = p\Gmax_m = 0$. The resulting problem \eqref{eqn:problemOrig} is a MI nonconvex problem which is not solvable in polynomial time. 
}

As mentioned in Section \ref{sec:intro}, applying BnB- or BnC-based solvers for \text{MI-SOCP} or \text{MI-SDP} problems to solve problem \eqref{eqn:problemOrig} has several drawbacks including a prohibitively high computational complexity for real-time applications. The following two approaches exhibit a significantly smaller computational burden.

\section{State-of-the-art two-step approach}\label{sec:2step}

In the two-step approach, that has been used by \cite{Marten} and \cite{Mamandur}, the integer constraints in \eqref{subeqn:stepConstr1} are continuously relaxed to turn problem \eqref{eqn:problemOrig} into the following nonconvex problem without discrete variables:
\begin{subequations}\label{eqn:problemContRelax}
\begin{eqnarray}
  & & (\ref{subeqn:objective})  \\
  & & \text{subject to} \nonumber \\  
  & & (\ref{subeqn:1stLFconstr}) - (\ref{subeqn:lastContConstr}) \\
  & & p\Gstep_{k,1} \leq p\Gstep_k \leq p\Gstep_{k,L_k}, \; \forall k \in \kK \label{subeqn:stepConstr1relax}
\end{eqnarray}
\end{subequations}
Although problem \eqref{eqn:problemContRelax} remains nonconvex, it represents a typical AC-OPF problem for which approximate solutions can be efficiently obtained using solvers such as the open-source IPOPT (see \cite{Molzahn} and references therein).

\begin{figure}[!ht]
\centering
\includegraphics[width=\columnwidth]{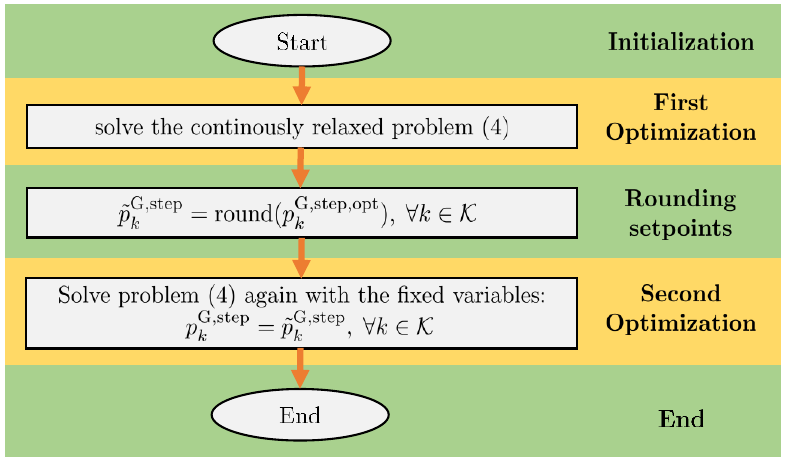}
\caption{Flowchart of the state-of-the-art two-step algorithm}
\label{fig:flowChart2step}
\end{figure}

In the first step of the two-step algorithm, problem \eqref{eqn:problemContRelax} is approximately solved yielding the solution of the continuously relaxed discrete setpoints $\bop\GstepOpt$. These relaxed setpoints are then rounded using a predefined problem-specific rule (e.g., upwards or downwards to the closest discrete setpoint). In the second step, problem \eqref{eqn:problemContRelax} is solved again with the discrete variables fixed at the values obtained by rounding to re-optimize the continuous optimization variables. The overall procedure is summarised in Fig.~\ref{fig:flowChart2step} where the $\text{round}(\cdot)$ operator rounds to the nearest discrete setpoint rather than to the nearest integer value.
\vspace{-3mm}
\vspace{-1mm}
\section{Proposed deflation approach}\label{sec:deflat}

To improve the results of the two-step approach, we propose the following iterative deflation approach. We first introduce the auxiliary binary variables $a_{k, l_k}$ to obtain the following equivalent reformulation of problem \eqref{eqn:problemOrig}:
\begin{subequations}\label{eqn:problemOrigDeflat}
%\begin{align}
%\begin{center}
\begin{eqnarray}
  & & \min_{\boldx, \boq\Gstep, \boa} f_2(\bop\G, \boa)  \label{subeqn:objectiveDeflat}\\
  & & \text{subject to} \nonumber \\
  & & p\G_m + \sum_{k \in \kK_m} \sum_{l_k \in \kL_k} a_{k, l_k} \cdot p\Gstep_{k, l_k} - p\D_m \nonumber \\
  & & = \sum_{n \in \kN_m} p_{mn}, \; \forall m \in \kN \label{subeqn:1stConstrDeflat} \\
  & & (\ref{subeqn:2ndLFconstr}) - (\ref{subeqn:lastLFconstr}) \label{subeqn:lastLFconstrDeflat} \\
  & & (\ref{subeqn:firstIneqConstr}) - (\ref{subeqn:lastContConstr}) \label{subeqn:lastContConstrDeflat} \\
  & & a_{k,l_k} \in \{0, 1\}, \; \forall l_k \in \kL_k, \; \forall k \in \kK, \label{subeqn:binConstr1} \\
  & & \sum_{l_k \in \kL_k} a_{k,l_k} = 1, \; \forall k \in \kK, \label{subeqn:binConstr2}
%\end{align}
\end{eqnarray}
%\end{center}
\end{subequations}
where the objective function $f_2(\cdot)$ is defined as follows
\begin{eqnarray}\label{eqn:objectiveDeflat}
  & & f_2(\bop\G, \boa) = \sum_{m \in \kN} \left| p\Gprior_m - p\G_m \right| \nonumber \\
%  & & \min_{\boldx, \boa} \; \sum_{m \in \kN}|p\Gprior_m - p\G_m| \nonumber \\
  % & & + \sum_{k \in \kK\G}|p\Gstepprior_k - \boa_k^T\bop\Gstep_k| \nonumber \\
  & & + \sum_{k \in \kK} \left |p\Gstepprior_k -  \sum_{l_k \in \kL_k} a_{k, l_k} \cdot p\Gstep_{k, l_k} \right|
\end{eqnarray}
and
\begin{subequations}
\begin{eqnarray}
    \boa_k & = & [a_{k,1}, \ldots, a_{k, L_k}]^T \\
    \boa & = & [\boa_1^T, \ldots, \boa_K^T]^T. \label{subeqn:boa}
\end{eqnarray}
\end{subequations}
The binary constraints in \eqref{subeqn:binConstr1} are continuously relaxed to turn problem \eqref{eqn:problemOrigDeflat} into the following nonconvex problem without the binary variables defined by \eqref{subeqn:binConstr1}:
\begin{subequations}\label{eqn:problemContRelaxDeflat}
\begin{eqnarray}
  & & \redtext{\eqref{subeqn:objectiveDeflat}} \\
  & & \text{subject to} \nonumber \\  
  & & (\ref{subeqn:1stConstrDeflat}) - (\ref{subeqn:lastContConstrDeflat}) \\
  & & 0 \leq a_{k,l_k} \leq 1, \; \forall l_k \in \kL_k, \; \forall k \in \kK, \label{subeqn:contRelaxConstr1} \\
  & & (\ref{subeqn:binConstr2}).
\end{eqnarray}
\end{subequations}
Note that also problems \eqref{eqn:problemContRelax} and \eqref{eqn:problemContRelaxDeflat} are equivalent.

\begin{figure}[!ht]
\centering
\includegraphics[width=\columnwidth]{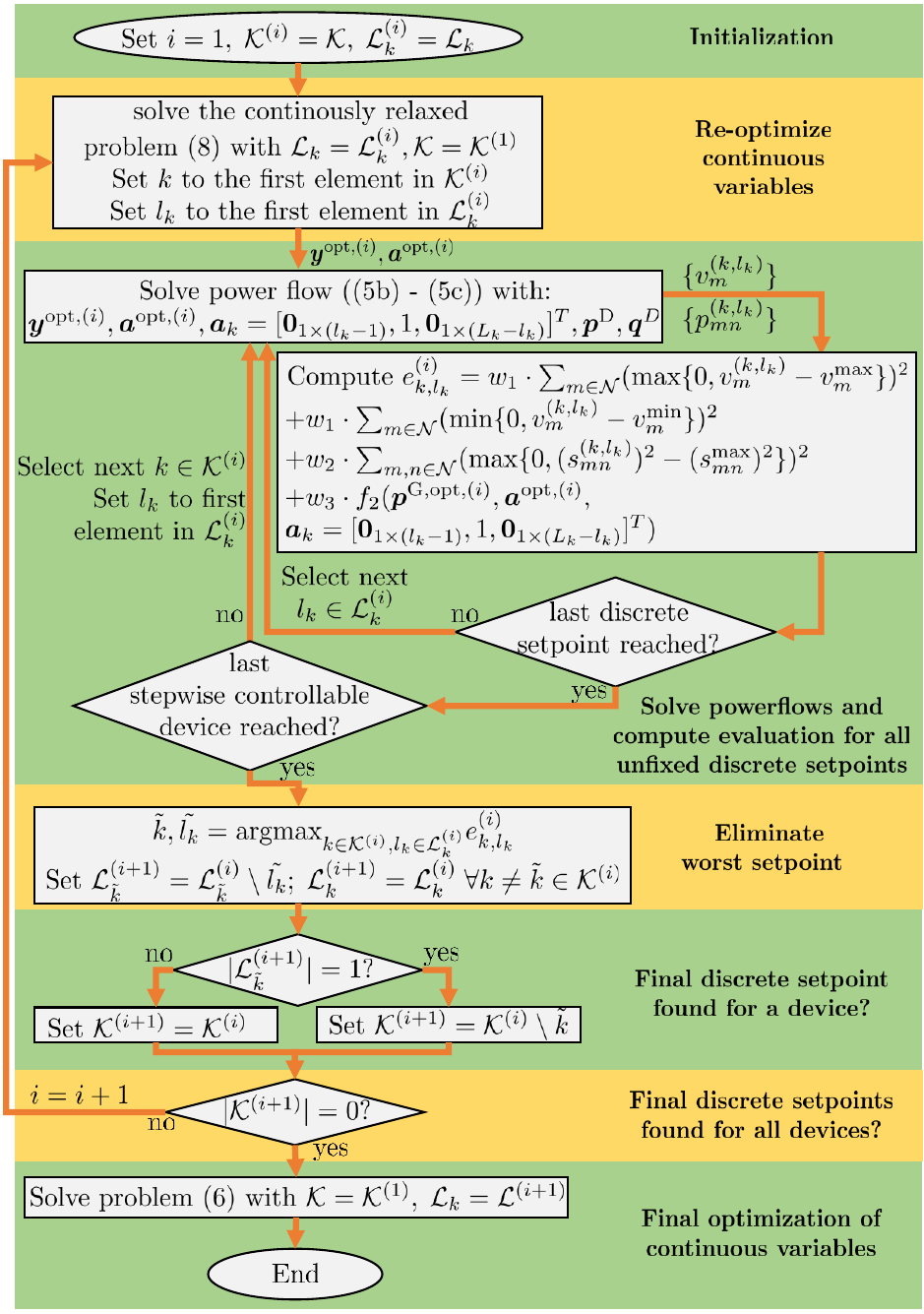}
\caption{Flowchart of the proposed deflation algorithm}
\label{fig:flowChartDeflation}
\end{figure}

The proposed iterative deflation procedure is illustrated in Fig.~\ref{fig:flowChartDeflation}. The first iteration is equivalent to the first step of the two-step procedure, where problem \eqref{eqn:problemContRelaxDeflat} is solved using a solver such as IPOPT. However, the fixation of discrete setpoints is performed more carefully in the deflation method. In each iteration, only one discrete setpoint is removed from a set of candidate discrete setpoints considered in problem \eqref{eqn:problemContRelaxDeflat}. After removal of a discrete setpoint, the next iteration starts with solving problem \eqref{eqn:problemContRelaxDeflat} again, now with the reduced set of candidate discrete setpoints. This procedure is repeated until only one discrete setpoint remains for each stepwise controllable device.
The selection of which candidate discrete setpoint to remove is guided by the proposed evaluation procedure, outlined in the second green box of Fig.~\ref{fig:flowChartDeflation}. For all candidate discrete setpoints remaining in the current iteration, and with the continuous setpoints obtained by solving problem \eqref{eqn:problemContRelaxDeflat} in the current iteration, a power flow is solved. The power flow result is used to compute the evaluation value $e_{k, l_k}$, which increases with the severity of voltage band and line loading constraint violations, as well as with the objective function value. The constraint violation and objective function value evaluations are combined using the weighting coefficients $w_1,$ $w_2$, and $w_3$. The candidate discrete setpoint with the largest, and thus worst, evaluation value is removed from the set of candidate discrete setpoints.

Note that the deflation method can be viewed as an extension to the two-step method rather than a replacement. In cases where the difference in the objective function value between the first and the second step of the two-step method is only marginal, there is only few potential for the deflation method to improve the result of the two-step approach. In this case, it is not worth spending additional computational load for performing the deflation method. When the above mentioned difference is large, however, the deflation method can be used to most likely obtain a significantly better solution.
\redtext{Similarly, when the two-step approach is not able to find a feasible solution solving the congestion, it is worth spending the additional computational effort of the deflation method to potentially obtain a feasible solution and solve the congestion.}

\redtext{In addition, note that we consider the stepwise controllability of legacy devices, such as old wind parks in Germany (see \cite{Gross}), only as an example of application and demonstration. On the HV level, DSOs often encounter other stepwise controllable devices, such as, e.g., switched capacitors or transformer tap changers, which are commonly deployed. The proposed deflation method can be readily adapted to solve MI-OPF problems involving these or even other stepwise controllable equipment. Thus, the method remains a highly relevant technique for DSOs when legacy devices have been replaced.}

% To Do: Streichkandidat:
%As shown in \cite{Marten}, transformer tap positions and DG feed-in can be jointly optimized using the two-step approach with satisfactory results. The reason why the two-step approach works well in this context can be attributed to the relatively small increments between transformer tap positions, typically around 1.75\,\%.
% To Do: 1,75%: How general is it?
%Hence, simply rounding to the closest discrete value is accurate enough, i.e., the discrete value is not fixed too early. In contrast, legacy devices and compensation devices (e.g., switched capacitors) often exhibit only very few steps within which they are controllable. In such cases, the proposed deflation scheme is expected to deliver more accurate results than the two-step approach due to its ability to handle the coarser step sizes more effectively.
%\redtext{Moreover, the proposed deflation method offers a significant potential for accuracy improvements when facing several stepwise controllable devices in the grid, rather than just one transformer, since the combinatorial aspect of the problem is more significant in this case.}

\section{Simulation results}\label{sec:simulations} 

% \begin{figure}[!ht]
% \centering
% \includegraphics[width=\columnwidth]{grid_new_hv.png}
% \caption{SimBench HV benchmark grid}
% \label{fig:grid}
% \end{figure}

\redtext{
We conduct our simulations in Python using the open-source tool pandapower \cite{Thurner} with the SimBench \cite{Meinecke} HV benchmark grid \verb|1-HV-urban--1-no_sw|.
% shown in Fig.~\ref{fig:grid}.
The benchmark grid operates at 110 kV and consists of 100 buses, to which 79 loads and 101 generators are connected. 
%The added lines each have a length of 20\,m in accordance with the line lengths of the already existing lines.
Among these, all 22 wind farms are considered controllable generators in our simulations.
The remaining generators, and all loads, are considered uncontrollable. The number of stepwise controllable generators among the controllable ones varies throughout the simulations. All remaining controllable generators are considered to be continuously controllable.
The storages contained in the original SimBench dataset are removed before the optimization since the optimization of storages across time steps is not the subject of this paper.
}

For the two-step method, we consistently rounded down to the nearest discrete setpoint. This specific rounding operator was chosen as it provided the best simulation results for the two-step approach, whereas alternative rounding strategies sometimes failed to yield feasible solutions.
The following weighting coefficients have been chosen for calculating the evaluation value $e_{k, l_k}$: $w_1 = 10^10 \frac{1}{\text{V}^2}$, $w_2 = 10^6 \frac{1}{\text{kW}^2}$ and $w_3 = 10 \frac{1}{\text{kW}}$. Hence, constraint violations are penalized more heavily to avoid infeasibility. The weights are chosen such that violations of line loading and voltage band constraints are equally penalized. The optimization problems are modeled using the open-source Python package Pyomo \cite{Bynum}. To solve the AC-OPF problems, we employ the open-source solver IPOPT (see [6] and references therein). The simulations are performed on a PC equipped with an Intel i7 1165G7 @2.80\,GHz featuring 8\,CPUs and 16\,GB RAM. 

\begin{figure}[!ht]
\centering
\includegraphics[width=1.1\columnwidth]{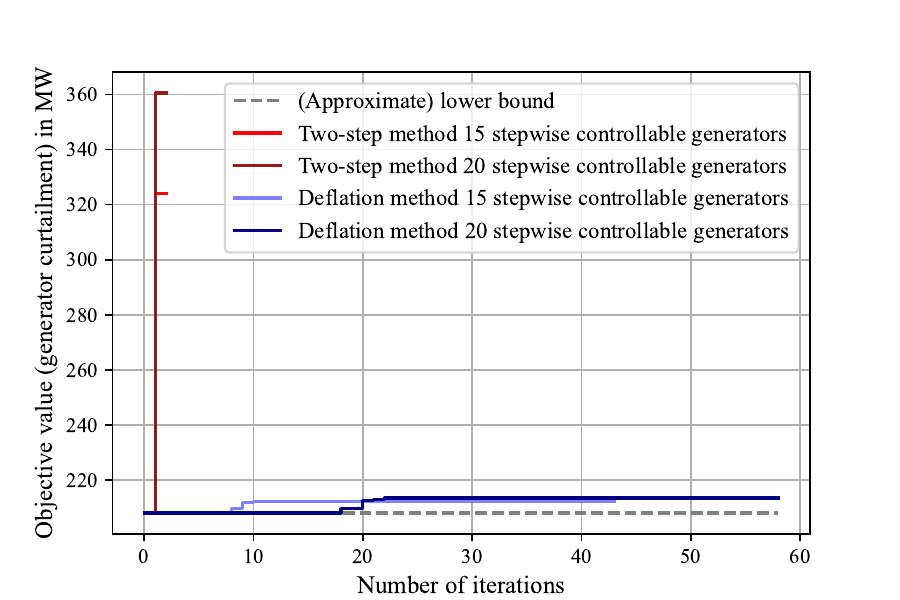}
\caption{Evolution of objective function value}
\label{fig:objFuncValIter}
\end{figure}

\redtext{
Fig.~\ref{fig:objFuncValIter} shows how the objective function value changes across iterations for both algorithms considering different numbers of stepwise controllable generators. Obviously, both algorithms start with the same objective function value which is unrealistically low since it is obtained by solving the continuously relaxed problem and is thus physically not attainable. In the second step of the two-step algorithm, all discrete setpoints are fixed at once leading to a large increase in the objective function value after re-optimization of the continuous variables. In the majority of the iterations of the deflation algorithm, the objective function value remains unchanged. This is due to the fact that often, when binary variables are eliminated, the exact same solution of the continuously relaxed discrete setpoints can be obtained with a new combination of the remaining continuously relaxed binary variables.
Future work could focus on improving computational efficiency by eliminating the binary variables in these iterations without re-solving the power flows and the continuously relaxed problem (\ref{eqn:problemContRelaxDeflat}).
During the remaining iterations, the objective function value gradually increases with each iteration as problem \eqref{eqn:problemContRelaxDeflat} is tightened from iteration to iteration. Once all discrete setpoints are fixed, the final objective function values are significantly lower than those obtained with the two-step procedure.
}

For the remainder of the simulations section, the number of stepwise controllable generators is increased from 1 to 22, which corresponds to decreasing the number of continuously controllable generators from 21 to 0. It is assumed that the generators can be curtailed to 0\,\%, 30\,\%, and 60\,\% (and 100\,\%) of their nominal power. All results are averaged over 20 time steps of the SimBench time series corresponding to the selected SimBench grid. Only time steps where voltage- and/or current-related congestions actually occurred were included in the analysis.
%To ensure congestions to occur, the maximum allowed line loading has been set to 25\,\%. This adjustment was necessary because, under the given SimBench grid and time series, congestion did not arise. In real-world scenarios, however, the rapid increase in DGs connected to HV grids is expected to cause a substantial rise in congestion incidents in the near future. 

\begin{figure}[!ht]
\centering
\includegraphics[width=1.1\columnwidth]{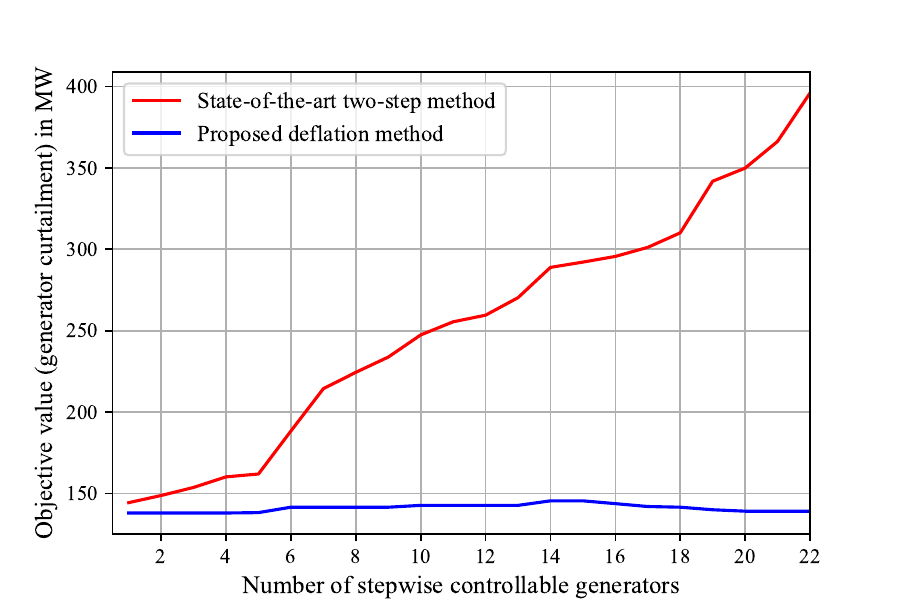}
\caption{Power curtailment vs. number of stepwise controllable generators}
\label{fig:objFuncVal}
\end{figure}
%\vspace{-1em}
% To Do: Wenn wir mehr Platz brauchen, könnten wir die Boxplot-Grafik + Beschreibung löschen.
\begin{figure}[!ht]
\centering
\includegraphics[width=1.1\columnwidth]{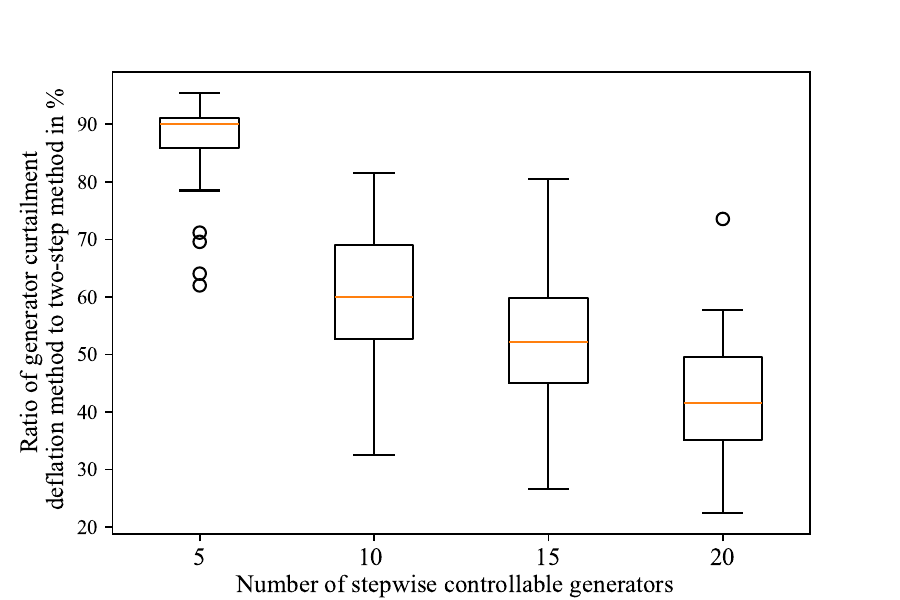}
\caption{Visualisation of improvement achieved by deflation method}
\label{fig:boxplot}
\end{figure}

Fig.~\ref{fig:objFuncVal} displays the obtained objective function value for both methods as a function of the number of stepwise controllable generators. It can be observed that the deflation method yields a significantly reduced required power curtailment compared to the two-step method. The improvement increases with the number of stepwise controllable generators. For a better comparability, the objective value obtained by the deflation method is normalized by dividing it by the value obtained from the two-step approach, representing the relative reduction in curtailed power. The results are shown in Fig.~\ref{fig:boxplot} using boxplots. As shown in the figure, the proposed deflation method reduces the curtailed power to as little as 30\,\% of the amount required by the two-step method. The system is operated closer to its operational limits, without violating them, and therefore uses the available capacities more efficiently within the operational boundaries of the grid.

\begin{figure}[!ht]
\centering
\includegraphics[width=1.0\columnwidth]{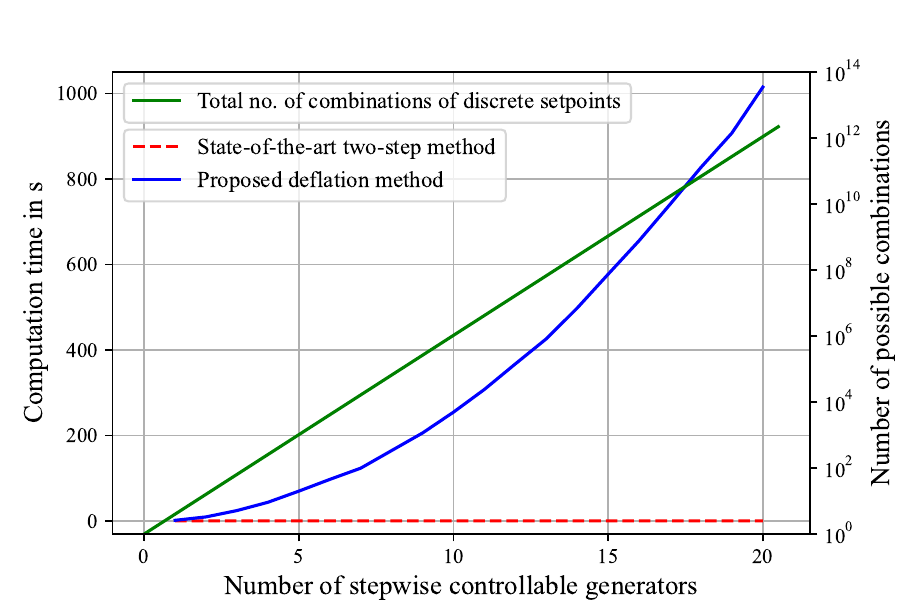}
\caption{Computation time vs. number of stepwise controllable generators}
\label{fig:calctime}
\end{figure}

% To Do: Diese Beschreibung können wir bei Bedarf kürzen:
Fig.~\ref{fig:calctime} shows the calculation time for the two approaches (left-hand side y-axis) and the number of possible combinations of discrete setpoints (right-hand side y-axis) as a function of the number of stepwise controllable generators. The calculation time for the two-step method remains constant since there are always two optimization problems to solve. For the deflation method, the number of instances of problem~\eqref{eqn:problemContRelaxDeflat} to be solved increases linearly with the number of candidate discrete setpoints since the number of $a_{k,l_k}$ also increases linearly, cf. Fig.~\ref{fig:flowChartDeflation}. However, the number of required power flow calculations grows quadratically with the number of candidate discrete setpoints. Although solving the optimization problem dominates the calculation time, the contribution from power flow computations is not negligible. Hence, the calculation time of the deflation method does not increase purely linearly but includes a minor quadratic component. Nevertheless, comparing the calculation time with the number of possible combinations (blue curve) reveals that this quadratically growing part is not crucial and the calculation time can be kept relatively low.

Note that the high computation time of around 16 minutes (1.000 seconds) is associated with 20 stepwise controllable generators not 20 controllable generators in total. In total there are 50 generators of which, in this case, 30 are continuously controllable. The number of 20 stepwise controllable generators is relatively high and usually does not occur in real grids. Thus, the proposed deflation method remains suitable for real time applications.

% \begin{center}
%     \begin{tabular}{ |c|c| } 
%     \hline
%         number of loads & calculation time deflation ( \\
%     \hline
%         1 discrete load & 12s \\
%         2 discrete loads & 28s \\   
%         4 discrete loads & 60s \\ 
%         8 discrete loads & 180s \\ 
%         16 discrete loads & 633s \\ 
% \hline
% \end{tabular}

%     \begin{tabular}{ |c|c| } 
%     \hline
%         number of loads & calculation time two-step  \\
%     \hline
%         1 discrete load & 6s \\
%         2 discrete loads & 3.4s \\   
%         4 discrete loads & 2.7s \\ 
%         8 discrete loads & 2.3s \\ 
%         16 discrete loads & 1.8s \\ 
% \hline
% \end{tabular}

% \end{center}

\section{Conclusion and outlook}\label{sec:conclusion}

\redtext{
In this paper, we have considered an AC-OPF problem with stepwise controllable generators. The discrete setpoints turn the nonconvex AC-OPF problem into a nonconvex MI-AC-OPF problem which cannot be solved optimally in affordable time. We have proposed an efficient iterative deflation procedure to obtain an approximate solution at a relatively low computational complexity. In each iteration of the deflation algorithm, one candidate discrete setpoint is removed from the set of all possible discrete setpoints, and a continuously relaxed version of the MI-AC-OPF problem is approximately solved afterward to re-optimize the remaining variables. The removed candidate discrete setpoint is the one yielding the worst power-flow-based evaluation. The computational complexity of the proposed algorithm grows only linearly with the number of discrete setpoints. Our simulation results demonstrate that the results of the state-of-the-art two-step approach can be substantially improved using the proposed method.
}

\redtext{
In future work, the performance of the proposed deflation method will be benchmarked against a commercial MINLP solver such as, e.g., Knitro. In addition, alternative evaluation approaches can be investigated to improve the selection process for determining which candidate setpoint to remove next. It can potentially be mathematically proven that in a large amount of iterations, candidate setpoints can be eliminated without re-solving the power flow and the continuously relaxed MI-AC-OPF problem. This will be the subject of future work and has the potential to drastically reduce the computational load of the proposed deflation method.
Also, a convex relaxation of the AC-OPF problem could be employed, enabling the relaxed MI-AC-OPF problem to be solved using solvers for convex MI programming problems. This would provide a lower bound for the original problem, offering insights into the optimality gap of the deflation method. The results of the deflation method could also be compared to those of MISDP- and MISOCP-based state-of-the-art approaches in terms of computational complexity, feasibility and accuracy. Finally, the proposed approach can be applied to other MI-OPF problems in power systems.
}

% \section*{Acknowledgment}

% The preferred spelling of the word ``acknowledgment'' in America is without 
% an ``e'' after the ``g''. Avoid the stilted expression ``one of us (R. B. 
% G.) thanks $\ldots$''. Instead, try ``R. B. G. thanks$\ldots$''. Put sponsor 
% acknowledgments in the unnumbered footnote on the first page.

\end{document}